# On Caccetta-Haggkvist Conjecture

Dhananjay P. Mehendale
Sir Parashurambhau College, Tilak Road, Pune-411030,
India

## Abstract

We show that we cannot avoid the existence of at least one directed circuit of length $\leq \left\lceil \dfrac{n}{r} \right\rceil$ in any digraph on $n$ vertices with out-degree $\geq r$. This is well-known Caccetta-Haggkvist problem.

**1. Introduction:** A digraph is a graph with directed edges, $G$, with vertex set $V(G) = \{v_1, v_2, \cdots, v_n\}$ and set of directed edges, $E(G)$. It is called simple when there does not exist any parallel directed edges or self loops, i.e. joining two vertices there is only one directed edge either from vertex $v_i$ to vertex $v_j$ or from vertex $v_j$ to vertex $v_i$, and, there is no directed edge from vertex $v_i$ to itself, for all $v_i, v_j \in V(G)$.

The number of directed edges emerging from a vertex $v_i$ is called the out-degree of $v_i$. The symbol $\lceil \ \rceil$ denotes the ceiling function. It is defined for argument $x$ as, $\lceil x \rceil$ = the smallest integer $\geq x$.

Caccetta-Haggkvist conjecture states that the existence of at least one directed circuit of length $\leq \left\lceil \dfrac{n}{r} \right\rceil$ in any digraph on $n$ vertices with out-degree $\geq r$ cannot be avoided. Many special cases of Caccetta-Haggkvist conjecture have been settled. It is settled for the out-degrees, $r = 2, 3, 4, 5$ and it is settled for out-degrees, $r \leq \sqrt{\dfrac{n}{2}}$. The case $r = 2$ was settled by Caccetta-Haggkvist, [1]. The case $r = 3$ by Hamidoune [2]. The case $r = 4$ and $5$ by Hoang and Reed [3]. The case $r \leq \sqrt{\dfrac{n}{2}}$ by Shen [4]. Shen's result [4] actually shows that for any $r$ the counter examples to the conjecture, if exist, are **finite**, thus, Caccetta-Haggkvist conjecture is **almost true!**

**2. Caccetta-Haggkvist Problem:** In this section we propose an algorithm which consists of labeling of vertices when out-degree condition is fulfilled. This algorithm will make it clear that because of out-degree condition, namely, the out-degree $\geq r$, for the digraph on $n$ points makes the appearance of a directed circuit of length $\leq \left\lceil \frac{n}{r} \right\rceil$ inevitable when the labeling process is continued up to last vertex. This is Caccetta-Haggkvist problem. To make the things clear and transparent we begin with an

**Example**: The case: $n = 7$ and $r = 3$. It will be evident after the discussion of this special case that the general case is just writing the special case **in terms of symbolic language!**

1) Take an unlabelled vertex and draw three directed edges emerging from it to three new unlabelled vertices and now assign label 1 to the vertex with out-degree three (Note that we label a vertex by next label only when three (since $r = 3$) outgoing edges will be emerging from it to three vertices, reaching to unlabeled vertices when they are available and to labeled ones when unlabeled vertices are not available), thus we have so far formed a figure containing in all four vertices, one vertex with out-degree = 3 and so labeled and other three unlabeled vertices with out-degree = 0, as shown in the following Fig. 1:

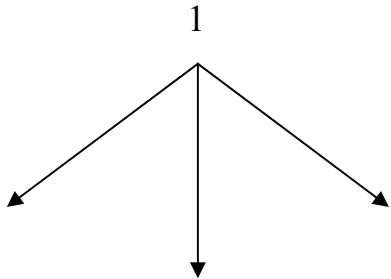

Fig. 1

2) Since $n = 7$, we have scope to add only three new vertices. So, we draw three directed edges emerging from an end vertex in the Fig. 2 to three new vertices and label thus formed vertex (with out-degree = 3) by next label 2, as shown in the following Fig. 2:

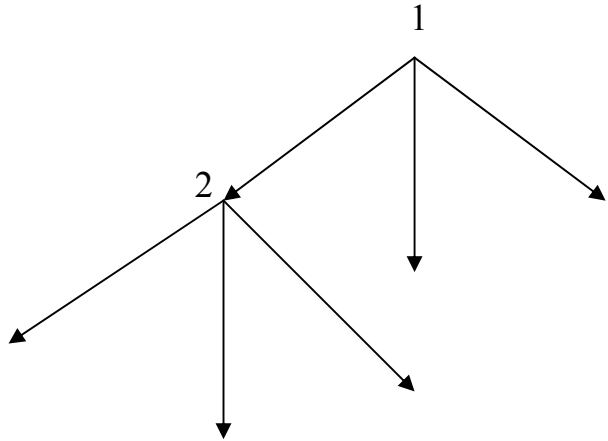

Fig. 2

3) Now we have in all 7 vertices in the digraph shown in the above Fig. 2. In this figure two vertices have got the labels, respectively, 1 and 2. Now, to achieve the out-degree condition we have to draw three outgoing edges from each so far unlabeled vertex reaching to some vertex among the 7 vertices. So far three 2-paths emerging from vertex with label 1 and passing through label 2 and ending in unlabeled vertex are formed and two 1-paths emerging from vertex 1 ending in an unlabeled vertex are formed. From an end vertex of a 2-path we have to take out three emerging edges. If we draw a directed edge among these edges connecting to any labeled vertex on the path then a closed circuit (a digon or triangle) will be formed. So, (to avoid this) we draw these edges ending into unlabeled vertices and label this vertex with the next label 3 as per our procedure as shown in the following Fig. 3:

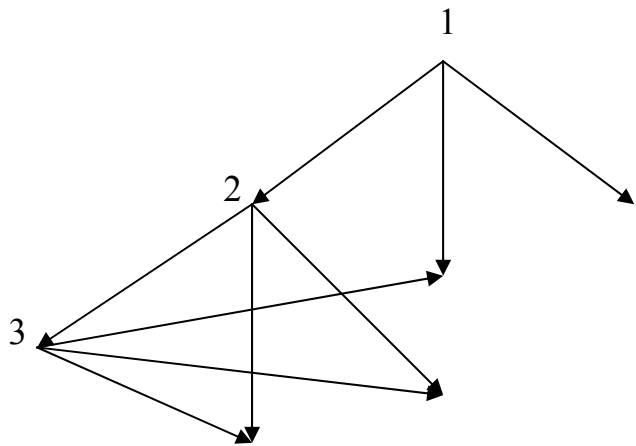

Fig. 3

4) We draw edges from next unlabelled vertex to unlabeled vertices as far as possible as shown in the following Fig. 4:

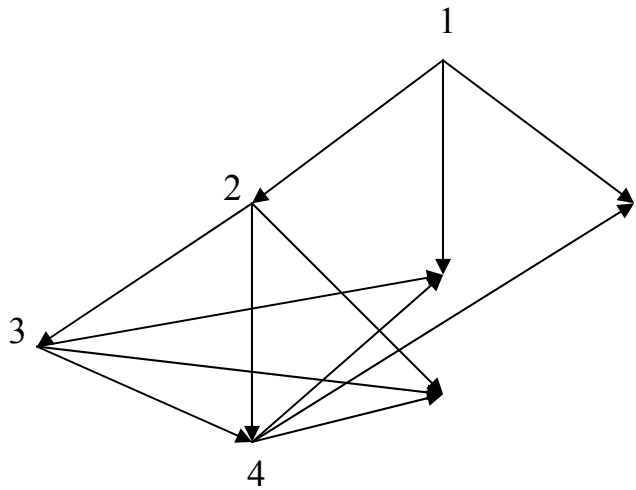

Fig. 4

5) To draw three edges from next unlabelled vertex to only unlabeled vertices is not possible since only two unlabelled vertices for entering of these new edges are now left. **So, for this case (and hereafter) we need to connect to some labeled vertex to fulfill the out-degree condition** and so the formation of a digon or triangle becomes inevitable as shown, as one of the inevitable choices as an example (a triangle in the colored lines) in the following Fig. 5:

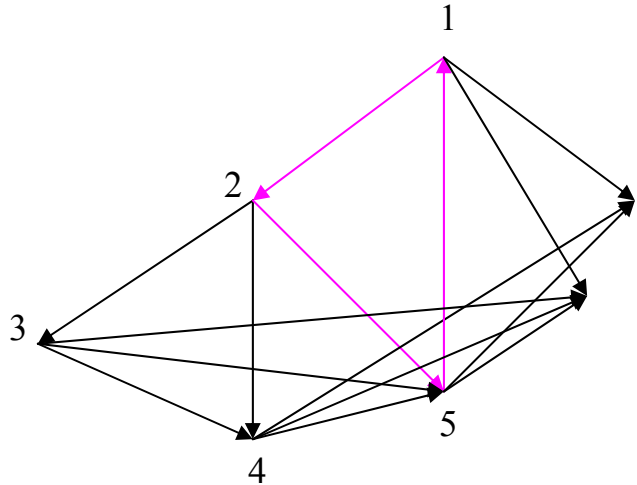

Fig. 5

**Remark 2.1:** There can be more than one unlabeled vertices to choose to assign a label by emerging out the number of directed edges equal to out-degree from them. The availability of multiple choices leads to different possible digraphs that get formed at the end of the process when each vertex gets a label.

One more **Example**: we now build a digraph with $n = 18$ and $r = 3$ by certain choice for labeling of unlabeled vertices:

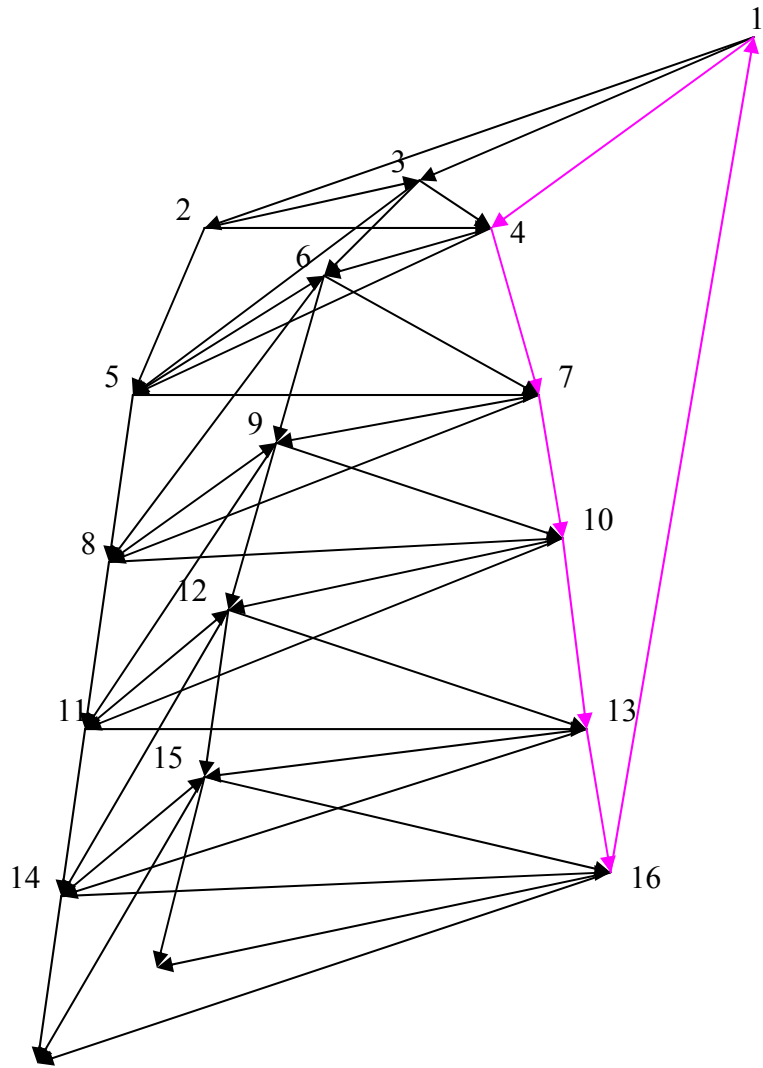

Fig. 6

It is clear to see that to label a vertex as 16 we need to extend a directed edge towards some already labeled vertex. **Longest among the shortest paths** from labeled vertices to the vertex with label 16 is available from vertex with label 1 so by joining the edge from vertex to be labeled as 16 to the vertex labeled 1 we see that we can't avoid the circuit of length $\left\lceil \dfrac{n}{r} \right\rceil = 6$, as shown by colored lines.

**Remark 2.2:** This algorithm of labeling of vertices in sequence (each time only after out-degree condition is fulfilled by the chosen vertex under consideration) by joining by directed edges only to unlabeled vertices cannot be continued indefinitely and a stage (i.e. a vertex to be labeled by next label in sequence) will always arrive for $n$ points and $r$ out-degree digraphs where for the left out vertices to be labeled after that stage we will need to create, in total, **directed edges joining** to at least $\left(\dfrac{r(r+1)}{2}\right)$ (already) **labeled vertices**. In this process of creating labeled vertices it is clear by the nature of the process that at each labeling of a vertex there appear $r$ new directed edges, creating $r$ new directed paths, bigger in length by unit, from each path reaching that vertex from the earlier labeled vertices. Also, it is easy to check that our procedure implies the existence of **directed path from each labeled vertex labeled earlier to each vertex labeled later**.

We now proceed formally with the steps of the algorithm, called **"labeling algorithm"**. This algorithm implies the unavoidability of the appearance of a directed path of length $\left\lceil \dfrac{n}{r} \right\rceil - 1$, for the case of $n$ vertices and $r$ out-degree, when we reach the stage of labeling a vertex with label $(n-r+1)$ which requires using already labeled vertex for entering of a directed edge from the vertex to be labeled by label $(n-r+1)$ and this leads to appearance of a directed circuit of length $\leq \left\lceil \dfrac{n}{r} \right\rceil$.

**Labeling Algorithm:**

1) Take an unlabeled vertex.
2) Draw $r$ directed edges emerging from this unlabeled vertex to newly taken $r$ unlabeled vertices and then assign label 1 to this vertex taken in step 1). (Hereafter, we assign next label to a vertex after we will draw $r$ directed edges emerging from it reaching other vertices.)
3) Choose any vertex among the unlabeled vertices to assign next label to it by drawing directed edges emerging from this vertex to some (or all) other unlabeled vertices (as far as possible) and to some (or one) new unlabeled vertices (vertex) taken newly at this time. (Note that adding a "new" vertex to already existing set of taken vertices is only allowed till the total count of so far created vertices is $\leq n$.)

4) Continue step 3) of labeling vertices by next labels in succession without taking as far as possible already labeled vertices for joining by directed edges entering in them and proceed choosing unlabeled vertices as far as possible and choosing labeled vertices for joining when unlabeled vertices exhaust till every vertex (among the $n$ vertices taken in all) gets a label. □

**Lemma 2.1:** For digraphs on $n$ points and out-degree $= r$, labeled using labeling algorithm there will always exist a directed path of length at most equal to $\left\lceil \dfrac{n}{r} \right\rceil - 1$ from any already labeled vertex to the vertex to get the label $(n - r + 1)$ for which we will be compelled of using some already labeled vertex to be joined by a directed edge to fulfill the out-degree condition for this vertex to get the label $(n - r + 1)$.

**Proof:** If we will execute the algorithm for out-degree $r = 1$, we can see that the stage of using already labeled vertex for entering of the directed edge will arrive for vertex $(n - r + 1) = n$, and up to this stage the (only) directed path of length $\left\lceil \dfrac{n}{r} \right\rceil - 1 = (n - 1)$ will be formed.

If we will execute the algorithm for out-degree $r = 2$, we can see that the stage of using already labeled vertex for entering of the directed edge will arrive for vertex $(n - r + 1) = n - 1$. Since out-degree $= 2$ for every labeled vertex there will be two distinct labeled vertices with higher labels adjacent to it. If we will choose, among the adjacent vertices, vertices with larger labels to form a directed path we can see that the vertex with label 1 will form a 1-path, $1 \to 3$ (or higher label). The vertex with labels in this path will extend as $1 \to 3$ (or higher label) $\to 5$ (or higher label). Thus, the minimum separation between the labels of the adjacent vertices on this path will be at least equal to $r = 2$ and so the length of the directed path will be equal to $\left\lceil \dfrac{n}{r} \right\rceil - 1 = \left\lceil \dfrac{n}{2} \right\rceil - 1$.

If we will execute the algorithm for out-degree $r$ we can see that same argument holds good and the result will be the formation of a directed path of length equal to $\left\lceil \dfrac{n}{r} \right\rceil - 1$ when executing the labeling algorithm we will be reaching to a vertex to which label $(n - r + 1)$ is to be assigned (even by proceeding using the smallest possible labels to define adjacencies forming a path like: $1 \to (r+1) \to (2r+1) \to \cdots \to (n-r+1)$). □

**Remark 2.3:** To view given unlabeled digraph on $n$ points and out-degree equal to $r$ as arrived at by so called "labeling algorithm", in which using unlabeled vertices for entering of directed edges emerging from the vertex to get a next label is followed as long as possible in its construction, what is required is to find **longest directed path** in the given unlabeled digraph **which has no backward directed edges** emerging from the vertices in the path which arrive later and reaching the vertices in the path which arrive earlier as one moves along the path from vertex to vertex along the direction of directed edges. Once we discover this path we assign label 1 to the vertex from which this path emerges after noticing the $r$ outgoing edges emerging from it among which there will be the vertex which is second vertex on this longest path. We choose this second vertex to assign label 2 proceeding as per labeling algorithm. We continue on these lines till every vertex in the unlabeled graph gets label as per labeling algorithm.

**Theorem 2.1(Caccetta-Haggvist Conjecture):** Every digraph on $n$ vertices with out-degree $\geq r$ contains a directed circuit of length $\leq \left\lceil \dfrac{n}{r} \right\rceil$.

**Proof:** The simple fact, namely, the impossibility of continuing the process of labeling vertices by "labeling algorithm" to the end without taking directed edges entering in some already labeled vertex at the stage of assigning to a vertex the label $(n-r+1)$, and existence of a directed path of length $\leq \left\lceil \dfrac{n}{r} \right\rceil - 1$ reaching to vertex to get a label $(n-r+1)$ from any previously labeled vertex gives rise to directed circuit of length $\leq \left\lceil \dfrac{n}{r} \right\rceil$.
□

**Remark 2.4:** It is clear to see that the case of out-degree $\geq r$ instead of out-degree $= r$ brings the stage, where one is enforced of using already labeled vertices as entering vertices for emerging directed edges from a vertex to get a label, will come earlier.

**3. Forward and Backward Edges and Differences:** We define and use the idea of forward and backward edges. We show that this simple idea settles at once a special case of the conjecture.

Let the vertex set of digraph under consideration be $V(G) = \{v_1, v_2, \cdots, v_n\}$. Let $X_0 = \{1, 2, \cdots, n\}$ be the set of suffixes of this assign label 1 to the vertex from which this path emerges vertex set. Let $D = \{1, 2, \cdots, n-1\}$ be the set of distinct differences of entries in the set $X_0$. We now state some simple definitions:

**Definition 3.1:** If there is directed edge from vertex $v_i$ to vertex $v_j$ and suppose $i < j$ then the directed edge from vertex $v_i$ to vertex $v_j$ will be called forward edge and the difference corresponding to a forward edge (between two distinct entries $i$ and $j$) in the set $X_0$, namely, $|i - j|$ will be called forward difference.

**Definition 3.2:** If there is directed edge from vertex $v_k$ to vertex $v_l$ and suppose $k > l$ then the directed edge from vertex $v_k$ to vertex $v_l$ will be called backward edge and the difference corresponding to a backward edge (between two distinct entries $k$ and $l$) in the set $X_0$, namely, $|k - l|$ will be called backward difference.

**Remark 3.1:** It is easy to see that in a complete symmetric digraph on $n$ points there are in all $n(n-1)$ directed edges out of which exactly $\frac{n(n-1)}{2}$ directed edges are actually forward edges producing exactly $\frac{n(n-1)}{2}$ forward differences and exactly $\frac{n(n-1)}{2}$ directed edges which are actually backward edges producing exactly $\frac{n(n-1)}{2}$ backward differences.

**Remark 3.2:** A simple digraph on the other hand can have at most $\frac{n(n-1)}{2}$ directed edges out of which some or all can be forward edges while the remaining are backward edges or vice versa, producing the corresponding forward and backward differences.

**Remark 3.3:** A simple digraph containing maximum number of directed lines (edges) and no directed circuit is the one containing all the $\frac{n(n-1)}{2}$ forward (or backward) edges. It is interesting to note that this digraph has out-degree equal to **zero**.

**Remark 3.4:** If we denote by $ex(p, C_n)$ the maximum number of directed lines a simple digraph on $p$ points can have without containing any closed directed circuit, $C_n$, of length 3 to $p$, then it is easy to see that $ex(p, C_n) = \frac{p(p-1)}{2}$. Since, a digraph on $p$ points containing all but only forward (or backward) directed edges can't contain any directed circuit and the cardinality of such set of forward (or backward) directed edges is $\frac{p(p-1)}{2}$.

**Remark 3.5:** A set of only forward (or backward) directed edges cannot form any directed circuit.

**Remark 3.6:** A simple digraph on $n$ points having out-degree $d$ is made up of in all (i.e. forward + backward) $nd$ directed edges, and therefore, incorporates in all $nd$ (forward + backward) differences.

**Remark 3.7:** For a simple digraph values of various differences and the cardinalities of their occurrences are as per the following table:

| Value of the difference | Number of Occurrences |
|---|---|
| 1 | (n-1) |
| 2 | (n-2) |
| ⋮ | ⋮ |
| (n-1) | 1 |

**Lemma 3.1 (Special case of Caccetta-Haggkvist Conjecture):** Let $G$ be a simple regular digraph on $n$ points and out-degree $r$ such that $\frac{(n-1)}{r} = 2$ and $\left\lceil \frac{n}{r} \right\rceil = 3$ then there exists a directed circuit of length 3 in $G$.

**Proof:** For this digraph $\frac{n(n-1)}{2} = nr$. Thus, for this graph the total number of directed edges and the number of all possible (forward + backward) differences together that can exist in a simple digraph are equal in number. Such digraph can be of type satisfying the following adjacency relations:

(i) There are $r$ out-going edges from vertex $v_1$ to each vertex in the set $S_1$ containing given $r$ vertices, where $S_1 = \{v_2, v_3, \cdots, v_{r+1}\}$.
(ii) There are $r$ out-going edges from vertex $v_2$ to each vertex in the set $S_2$ containing given $r$ vertices, where $S_2 = \{v_3, v_4, \cdots, v_{r+2}\}$.
(iii) There are $r$ out-going edges from vertex $v_2$ to each vertex in the set $S_3$ containing given $r$ vertices, where $S_3 = \{v_4, \cdots, v_{r+3}\}$.

$$\vdots$$

(iv) There are $r$ out-going edges from vertex $v_r$ to each vertex in the set $S_r$ containing given $r$ vertices, where $S_r = \{v_{r+1}, v_4, \cdots, v_{2r+1}\}$.
(v) Continue till you reach vertex $v_n$ and take $r$ out-going edges from vertex $v_n$ to each vertex in the set $S_n$ containing given $r$ vertices, where $S_n = \{v_1, v_2, \cdots, v_r\}$.

Since $\frac{(n-1)}{2} = r$, therefore, $\frac{(n+1)}{2} = r+1$ and $n = 2r+1$.

Thus, for this graph we have the following directed circuit of length 3:

$$v_1 \to v_{\frac{n+1}{2}} \to v_n \to v_1.$$

Thus, all the simple digraphs of this type containing $\frac{n(n-1)}{2} = nr$ directed edges and their isomorphs will contain a directed 3-circuit. Other digraphs containing $\frac{n(n-1)}{2} = nr$ directed edges and suppose there is no representative directed edge representing each difference among total possible $\frac{n(n-1)}{2}$ distinct difference then there is repetition of some difference more than the possible occurrences for a simple digraph and thus such digraph will contain a digon. Hence, etc. □

To avoid the existence of directed circuits actually one should be able to select directed edges of some one type, namely, forward or backward directed edges. We see that it is impossible to make such choice while maintaining the constraint on the out-degree, like, the out-degree $\geq r$.

**A Method of Digraph Construction with Maximum Possible Forward Edges:** We discuss procedure to construct simple digraphs of out-degree $\geq r$ in most general way using only forward directed edges as far as possible, and each time taking care of avoiding self loops or parallel edges, as follows:

(i) We choose some $r$ out going edges from vertex $v_1$ to each vertex in some subset $T_1 \subseteq S_1$ containing some $\geq r$ vertices, where $S_1 = \{v_2, v_3, \cdots, v_n\}$.

(ii) We Choose some $r$ out going edges from vertex $v_2$ to each vertex in some subset $T_2 \subseteq S_2$ containing some $\geq r$ vertices, where $S_2 = \{v_3, \cdots, v_n\}$.

(iii) We choose some $r$ out going edges from vertex $v_3$ to each vertex in some subset $T_3 \subseteq S_3$ containing some $\geq r$ vertices, where $S_3 = \{v_4, \cdots, v_n\}$.

(iv) We continue till we reach vertex $v_n$ and form $S_n$, and further choose some $r$ out going edges from vertex $v_n$ to each vertex in some subset $T_n \subseteq S_n$ containing some $\geq r$ vertices, where $S_n = \{v_1, \cdots, v_{n-1}\}$.

It is easy to check that this procedure of choosing forward directed edges each time to avoid backward edges cannot be continued indefinitely, in fact, not beyond the vertex with label $v_{n-r}$. And in order to fulfill the out-degree $\geq r$ we are forced to choose a backward edge emerging from vertex with label $v_{n-r+1}$. Even if we restrict the out-degree to be exactly $= r$ and choose adjacencies defined through (i)-(v) on page 11, then even for this digraph we cannot avoid the directed circuit

$$v_1 \to v_{r+1} \to v_{2r+1} \to v_{3r+1} \to \cdots \to v_{kr+1} \to v_1, \text{ where } k = \left\lceil \frac{n}{r} \right\rceil.$$

Thus, a smaller circuit of length $k = \left\lceil \frac{n}{r} \right\rceil$ can't be avoided.

**4. Transparency Matrix for Digraphs:** We define the so called transparency matrix for digraphs as was done for graphs in [5], and obtain two proofs for Caccetta-Haggkvist conjecture. We also provide a simple inductive proof for Seymour's second neighborhood conjecture [6].

**Definition 4.1:** Transparency matrix, $T(G)$, associated with a digraph $G$ containing $p$ points and $q$ directed lines, is the following $p \times p$ matrix:

$$T(G) = [a_{ij}]_{p \times p}$$

where, $a_{ii} = 0$, $a_{ij} = k$, where $k$ is the distance between vertices $v_i$ and $v_j$, i.e. it is the length of the shortest path of directed edges starting from vertex $v_i$ and ending in vertex $v_j$.

**Remark 4.1:** When $a_{ij} = k$ it meant that there is a directed path which consists of directed edges, directed in the same direction, where the first directed edge emerges from vertex $v_i$ and the last or $k^{th}$ edge ends in the vertex $v_j$.

**Some Interesting Properties of $T(G)$:**

(1) When there is no directed path connecting vertices $v_i$ and $v_j$ then $a_{ij} = \infty$.
(2) When there is a directed edge from vertex $v_i$ to $v_j$ then $a_{ij} = 1$.
(3) If we replace all the so called distances $k \geq 2$ by zero then the transparency matrix becomes the usual adjacency matrix, $A(G)$ for the digraph.
(4) All the two by two principle sub-matrices of the form $\begin{bmatrix} 0 & 1 \\ 1 & 0 \end{bmatrix}$ form the directed double-edges or digons in $G$.

(5) Only those pairs of vertices $v_i, v_j$ are available for contraction (by putting vertex $v_i$ onto $v_j$) for which the principle sub-matrix formed by the elements in the intersection of $i$-th row/column and $j$-th row/column has the form $\begin{bmatrix} 0 & 1 \\ 0 & 0 \end{bmatrix}, \begin{bmatrix} 0 & 0 \\ 1 & 0 \end{bmatrix}, \begin{bmatrix} 0 & 1 \\ 1 & 0 \end{bmatrix}$.

(6) Only the diagonal elements of $T(G)$ are zero and all the other elements of $T(G)$ are greater than zero.

(7) In one contraction obtained by identifying some two adjacent vertices, i.e. vertices at unit distance, the size of $T(G)$ reduce by one unit.

We now look at the effect of a contraction of some directed edge on the transparency matrix. When there is a directed edge $v_i \to v_j$ then by $(v_i \Rightarrow v_j)$ we denote the contracting of edge $(v_i, v_j)$ and identifying the vertex $v_i$ with the vertex $v_j$. Let $\tilde{G}$ be the digraph that results after the operation $(v_i \Rightarrow v_j)$ on $G$, and let $T(\tilde{G})$ denotes the transparency matrix for $\tilde{G}$, then

$$T(\tilde{G}) = [a_{ij}]_{(p-1)\times(p-1)}$$

can be obtained from $T(G)$ by performing the following operations:

(1) Replace all elements $a_{jk}$, $k \neq j$ by min $\{a_{ik}, a_{jk}\}$.
(2) Replace all elements $a_{kj}$, $k \neq j$ by min $\{a_{ki}, a_{kj}\}$.
(3) Delete $i$-th row and $i$-th column.
(4) If the shortest directed path (deciding distance) joining vertices $v_m$ and $v_n$ contains the edge $(v_i, v_j)$ then replace the entry $a_{mn}$ by $(a_{mn} - 1)$ in all such $v_m$ and $v_n$.
(5) Keep all other elements as they are.

**Remark 4.2:** How one finds out whether the condition mentioned in (4) is true or not? The answer is simple: The condition will be true if and only if $a_{mn} = a_{mi} + a_{jn} + 1$.

**Remark 4.3:** It is easy to see that if $a_{ij} + a_{ji} = k$, where $(a_{ij}, a_{ji})$ forms complementary pair of elements of transparency matrix, then there will exist a directed circuit of length $\leq k$ in the digraph $G$.

**Definition 4.2:** The regular digraphs obtained by choosing the specific method of defining adjacencies as defined in steps (i)-(v) in the proof of lemma 3.1 will be called **uniform regular digraphs**.

**Definition 4.3:** The digraphs of out-degree $\geq r$ obtained by choosing the specific method of defining adjacencies as defined in steps (i)-(v) in the proof of lemma 3.1 with additional adjacencies chosen in exact succession, like $v_{r+2}, v_{r+3}, \cdots$ in the first set, $v_{r+3}, v_{r+4}, \cdots$ in the second set, etc. when the out-degree for that vertex $> r$, will be called **uniform digraphs of out-degree** $\geq r$.

**Definition 4.4:** The regular digraphs of out-degree $r$ obtained by any other choice of defining adjacencies other than defined in steps (i)-(v) in the proof of lemma 3.1 will be called **nonuniform regular digraphs**.

**Definition 4.5:** The digraphs of out-degree $\geq r$ obtained by any other choice of defining adjacencies other than defined in steps (i)-(v) in the proof of lemma 3.1 with additional adjacencies **not** chosen in exact succession, like $v_{r+2}, v_{r+3}, \cdots$ in the first set, $v_{r+3}, v_{r+4}, \cdots$ in the second set, etc. when the out-degree for that vertex $> r$, will be called **nonuniform digraphs of out-degree** $\geq r$.

For uniform regular digraphs Caccetta-Haggkvist conjecture is true and as mentioned on page 6 we can see the existence of many directed circuits, like the one mentioned below:

$$v_1 \to v_{r+1} \to v_{2r+1} \to v_{3r+1} \to \cdots \to v_{kr+1} \to v_1, \text{ where } k = \left\lceil \frac{n}{r} \right\rceil.$$

We are now ready to state the first version of the Caccetta-Haggkvist conjecture in the language of entries of the transparency matrix defined for digraphs as follows:

**Theorem 4.1 (Caccetta-Haggkvist Conjecture):** In a digraph on $n$ points and of out-degree $\geq r$ there exist at least one pair of matrix elements of the transparency matrix, $(a_{ij}, a_{ji})$ such that $a_{ij} + a_{ji} \leq \left\lceil \frac{n}{r} \right\rceil$.

**Lemma 4.1:** For regular uniform digraphs on $n$ points and of out-degree $r$
$$\min\{a_{ij} + a_{ji}\} = \left\lceil \frac{n}{r} \right\rceil,$$
where $a_{ij}, a_{ji} \in T(G)$, the associated transparency matrix.

**Proof:** The special method of defining adjacency relations fulfilled by the transparency matrix it is clear to see that in the uniform digraphs maximum number of directed circuits of length $\left\lceil \frac{n}{r} \right\rceil$ and are formed which are uniformly distributed over the digraph. The result is clear from the special form of the transparency matrix. □

**Lemma 4.2:** For regular nonuniform digraphs on $n$ points and of out-degree $r$ we have $\min\{a_{ij} + a_{ji}\} \leq \left\lceil \frac{n}{r} \right\rceil$, where $a_{ij}, a_{ji} \in T(G)$, the associated transparency matrix.

**Proof:** For the nonuniform regular digraphs of out-degree $r$ we randomly choose some $r$ vertices to emerge from each vertex in the vertex set to form all possible such digraphs. This essentially leads to **replacement** of one or more directed circuits of length $\left\lceil \frac{n}{r} \right\rceil$ by some directed circuits of smaller length leading to: $\min\{a_{ij} + a_{ji}\} \leq \left\lceil \frac{n}{r} \right\rceil$, □

**Lemma 4.3:** For uniform digraphs of out-degree $\geq r$ on $n$ points we have $\min\{a_{ij} + a_{ji}\} \leq \left\lceil \frac{n}{r} \right\rceil$, and difference, $|\min\{a_{ij} + a_{ji}\} - \max\{a_{ij} + a_{ji}\}|$, is minimum, where $a_{ij}, a_{ji} \in T(G)$, the associated transparency matrix.

**Proof:** Clear from the special form of the transparency matrix. □

**Lemma 4.4:** For nonuniform digraphs of out-degree $\geq r$ on $n$ points we have $\min\{a_{ij} + a_{ji}\} \leq \left\lceil \frac{n}{r} \right\rceil$, where $a_{ij}, a_{ji} \in T(G)$, the associated transparency matrix.

**Proof:** For the nonuniform digraphs of out-degree $\geq r$ we randomly choose some $\geq r$ vertices to emerge from each vertex in the vertex set to form all possible such digraphs. This essentially leads to **replacement** of one or more directed circuits of length $\approx \left\lceil \dfrac{n}{r} \right\rceil$ by some directed circuits of smaller length. Hence etc. □

**Proof of theorem 4.1:** Follows from the above lemmas. □

**Definition 4.6:** The number of vertices with value $n = 0, 1, \cdots, r$ will be called **first** range (since for all this range $\left\lceil \dfrac{n}{r} \right\rceil = 1$).

**Definition 4.7:** The number of vertices with value $n = r+1, r+2, \cdots, 2r$ will be called **second** range (since for all this range $\left\lceil \dfrac{n}{r} \right\rceil = 2$).

**Definition 4.8:** The number of vertices with value $n = (k-1)r+1, \cdots, kr$ will be called ***k*-th** range (since for all this range $\left\lceil \dfrac{n}{r} \right\rceil = k$).

**Proof of theorem 4.1 (Induction on the value of range):** We proceed by induction on the value of range of vertices.

**Step 1):** It is clear to see that the result holds for first, second, third range, i.e. for the value of range = 1, 2, 3.

**Step 2):** We assume the result by induction for the value of range = $k$, and proceed to prove it for the value of range = $(k+1)$.

Suppose the contrary, i.e. suppose the result is not true for the value of range = $(k+1)$. Therefore, there exist digraphs on $n = kr+1, \cdots, (k+1)r$ vertices for which $\min\{a_{ij} + a_{ji}\} > (k+1)$, where $a_{ij}, a_{ji} \in T(G)$, the associated transparency matrix for the digraph. Consider such a digraph on $n = kr+1$ points for which $\min\{a_{ij} + a_{ji}\} > (k+1)$. This implies that every directed circuit on this digraph is $> (k+1)$. If we perform one contraction the range due to this contraction will now reduce to $k$, and the length of some directed circuit will reduce by unit i.e. the length of directed circuit will be now $> k$, thus for $k$

□

range we will have $\min\{a_{ij} + a_{ji}\} > k$, a contradiction to induction hypothesis. Hence etc.

**Seymour's Second Neighborhood Conjecture:** We now proceed to give a simple inductive proof for this conjecture [6].

Seymour's Second Neighborhood Conjecture claims that for every oriented digraph there exists at least one vertex, $v$, such that $v$ has at least as many neighbors at distance two ($N_v^{++}$) as it has at distance one ($N_v^+$). We prove the following equivalence of this conjecture in the language of the associated transparency matrix with the oriented digraph under consideration:

**Theorem 4.2 (Matrix Version):** If we consider the transparency matrix, $T(G)$ associated with any oriented digraph on $n$ points, $G$ say, then the matrix will contain at least one row (column) such that the count of the matrix elements having value greater than or equal to two is grater than or equal to the count of the matrix elements having value exactly equal to one.

**Proof:** We proceed by induction on the number of vertices, $n$.

**Step 1:** For $n = 1, 2, 3$ the result is clear.

**Step 2:** We assume the result to be true for all number $\leq n$ and prove the same for the number ($n+1$).
        Suppose that the contrary is true. So, in every row associated with some vertex in the oriented digraph on ($n+1$) points the count of elements having value equal to two is less than the count of elements having value equal to one. Now choose any pair of adjacent vertices, $v_k, v_l$ at distance =1 and contract the edge ($v_k \to v_l$) by putting $v_k$ on $v_l$. In the language of transparency matrix this is equivalent to performing the contraction operation as described in 7) in the section where we discussed some interesting properties of $T(G)$. This leads to the formation of the transparency matrix associated with the resulted oriented graph on $n$ points.
        Now, note that due to this contraction operation we arrive at the following outcomes:
1) When $a_{ki}, a_{li} \geq 2$, there is no contribution towards increase in the count of elements having value greater than or equal to two in the resulted row (column) of the resulted transparency matrix after this contraction, i.e., (original count of elements with value greater than or equal to two) remains equal to (original count of elements with value greater than or equal to two).

2) When one of the $a_{ki}, a_{li} = 1$, there is no contribution towards increase in the count of elements having value greater than or equal to two in the resulted row (column) of the resulted transparency matrix after this contraction, i.e., (original count of elements with value greater than or equal to two) remains equal to (original count of elements with value greater than or equal to two). But on the contrary there is likelihood of increase in the count of elements having value equal to one, i.e., (original count of elements with value one) may changes to (original count of elements with value one +1).

3) When both i.e., $a_{ki} = 1$ and $a_{li} = 1$, there is no contribution towards increase in the count of elements having value greater than or equal to two in the resulted row (column) of the resulted transparency matrix after this contraction, i.e., (original count of elements with value greater than or equal to two) remains equal to (original count of elements with value greater than or equal to two). Also, there is no increase in the count of elements having value equal to one, i.e., (original count of elements with value one) remains (original count of elements with value one).

Thus, in any event there is no outcome which can cause increment in the count of elements having value greater than or equal to two, and thus the originally taken assumption, namely, in every row associated with some vertex in the oriented digraph on ($n+1$) points the count of elements having value greater than or equal to two is less than the count of elements having value equal to one, remains valid for this new oriented digraph on $n$ points arrived at due to one contraction remains valid if the original assumption is true. But, this contradicts the induction hypothesis. Hence this contrary supposition "in every row associated with some vertex in the oriented digraph on ($n+1$) points the count of elements having value equal to two is less than the count of elements having value equal to one" is not valid. Hence, etc.

5. **Cycles, Directed Circuits, and Cyclic Monomials:** Let the vertex set of digraph under consideration be $V(G) = \{v_1, v_2, \cdots, v_n\}$. Let $X_0 = \{1, 2, \cdots, n\}$ be the set of suffixes of the vertex set. The total combinations of symbols of $X_0$, taken $j$ at a time, are $\binom{n}{j}$ in number.

Let $X_1$ be the set of all possible (distinct) cycles of length $j$ made from some $j$ symbols of set $X_0$. We denote a cycle in this set of cycles, following the standard notation, as $(i_1, i_2, \cdots, i_j)$. Let $X_2$ be the set of all possible directed circuits in the digraph $G$ under consideration. We view a cycle $(i_1, i_2, \cdots, i_j)$ as a directed circuit $v_{i_1} \to v_{i_2} \to \cdots \to v_{i_j} \to v_{i_1}$ in the digraph $G$ under consideration. It is easy to see that the total number of such cycles of length $j$ made up of some $j$ symbols of $X_0$ will be $\binom{n}{j}(j-1)!$ in number. Therefore, the total number of all possible directed circuits of length $j$ in the digraph $G$ under consideration will be also equal to $\binom{n}{j}(j-1)!$ in number. Let $A(G) = [a_{ij}]_n$ be the adjacency matrix for digraph $G$. Let us associate a cyclic $j$-monomial, $a_{i_1 i_2} a_{i_2 i_3} \cdots a_{i_j i_1}$, made up of elements of $A(G)$, with each cycle $(i_1, i_2, \cdots, i_j)$ considered above. Let $X_3$ be the set of all possible cyclic $j$-monomials then, clearly, the total number of all possible cyclic $j$-monomials will be $\binom{n}{j}(j-1)!$ in number. Now, how many cyclic monomials are there such that all of them contain a non-diagonal element $a_{ij}$? It is easy to check that they are $\binom{n-2}{j-2}(j-3)!$ in number. Therefore, presuming the vanishing, with each setting to zero of some new non-diagonal element in $A(G)$, of new $\binom{n-2}{j-2}(j-3)!$ cyclic $j$-monomials (and not less due to repeated counting of already vanished monomials in the earlier setting to zero of some non-diagonal element of $A(G)$, the adjacency matrix) it is clear to see that as a rough estimate at least $\dfrac{n(n-1)(j-2)}{j}$ non-diagonal elements of $A(G)$ must be made zero to eliminate all the cyclic $j$-monomials. Since $A(G)$ of a simple digraph already contain half of its elements equal to zero, so the elements required to **make zero newly** are equal to $\dfrac{n(n-1)(j-2)}{2j}$ in number.

**Remark 5.1:** In actuality the required count of the elements of $A(G)$ to be made zero is **more** than $\dfrac{n(n-1)(j-2)}{j}$ and the eliminated $\binom{n-2}{j-2}(j-3)!$ monomials with each setting to zero of some non-diagonal element contains the already vanished cyclic monomials in the earlier setting to zero of some non-diagonal elements of $A(G)$.

**Remark 5.2:** When a non-diagonal element $a_{ij}$ is set to zero to eliminate certain number of $j$-monomials then setting $a_{ik}$ to zero such that $k \neq j$ ($a_{kj}$ to zero such that $k \neq i$) eliminates **same** number of cyclic $j$-monomials, since $a_{ij}$ and $a_{ik}$, ($a_{ij}$ and $a_{kj}$) cannot belong to same monomials, and so they eliminate different monomials.

**A New Version of Caccetta-Haggkvist Conjecture:** For the simplicity of presentation we take the following equivalent statement of the problem:

**Theorem 5.1:** Simple digraphs on $n$ points with minimum out-degree $\left\lceil \dfrac{n}{r} \right\rceil$ which do not contain directed circuits of length at most $r$ **do not exist**, i.e. simple digraphs on $n$ points which do not contain directed circuits of length at most $r$ have out-degree $< \left\lceil \dfrac{n}{r} \right\rceil$.

**Proof:** As per the above the above discussion the total number of directed circuits that can exist of length $r$ = total number of $r$-cycles = total number of cyclic $r$-monomials made up of elements of the adjacency matrix, $A(G)$, corresponding to the digraph on $n$ points. We have seen that this number $= \binom{n}{r}(r-1)!$. By setting any nondiagonal element, $a_{ij}$ to zero we eliminate $\binom{n-2}{r-2}(r-3)!$ cyclic monomials (i.e. their value becomes zero) which in effect eliminate the same number of directed circuits of length $r$ from the digraph. Thus, in order to eliminate all the directed circuits of length $r$ we should set to zero at least number of elements equal to ratio of the above two numbers, i.e. the ratio of total number of directed circuits, i.e. $\binom{n}{r}(r-1)!$, to the number of directed circuits of length $r$ that vanish due to setting some nondiagonal element, $a_{ij}$, to zero, i.e. $\binom{n-2}{r-2}(r-3)!$. It is easy to check that this ratio is actually equal to $\frac{n(n-1)(r-2)}{r}$. In a simple digraph when some element $a_{ij} = 1$, then $a_{ji} = 0$. Therefore, when $a_{ij} = 1$ the monomials containing its complement $a_{ji}$ are nonexistent because then $a_{ji} = 0$, already. Therefore, we actually need at least $\frac{n(n-1)(r-2)}{2r}$ elements of the adjacency matrix $A(G)$ to set to zero. Also, the adjacency matrix, $A(G)$, of a simple digraph can have at most $\frac{n(n-1)}{2}$ nonzero entries because of the same reason, namely when $a_{ij} = 1$, then $a_{ji} = 0$. Therefore, to make this digraph free of all the directed circuits of length $r$ its $A(G)$ can have at most $\frac{n(n-1)}{2} - \frac{n(n-1)(r-2)}{2r} = \frac{n(n-1)}{r}$ nonzero elements. Now, for this simple digraph to have out-degree at least equal to $\left\lceil \frac{n}{r} \right\rceil$ it must contain at least $n\left\lceil \frac{n}{r} \right\rceil$ directed edges, so, its $A(G)$ must contain at least $n\left\lceil \frac{n}{r} \right\rceil$ nonzero diagonal elements. We now proceed with the following easy claim which settles the theorem:

**Claim:** $n\left\lceil \dfrac{n}{r} \right\rceil > \dfrac{n(n-1)}{r}$

**Proof:** For $n = (k-1)r+1, (k-1)r+2, \cdots, kr$ we have $\left\lceil \dfrac{n}{r} \right\rceil = k$ and $\left(\dfrac{n}{r}\right) = (k-1)+\dfrac{1}{r}, (k-1)+\dfrac{2}{r}, \cdots, k$ and the claim and so the theorem is now clear. □

## Acknowledgements

Mr. Riko Winterle initiated me to this problem. I am extremely thankful to him for stimulating discussions on this problem through his emails. Many thanks are also due to Dr. M. R. Modak for useful discussions and encouragements.